%% file: Modifications.tex
\documentstyle{article}

\newtheorem{conjecture}{Conjecture}

\newtheorem{theorem}{Theorem}

\renewcommand{\refname}{\bf References}
\makeatletter
\renewcommand{\@evenhead}{\thepage \hfil {\leftmark} \hfil}
\renewcommand{\@oddhead}{\hfil {\rightmark} \hfil \thepage}
\makeatother
\input tcilatex


\begin{document}
\count0=1

\title{{\Large \textbf{Large Cycles in Graphs \\ Around Conjectures of Bondy and Jung - \\ Modifications and Sharpness}}}
\author{\normalsize  Zhora Nikoghosyan}

\maketitle

\begin{abstract}
A number of  new sufficient conditions for generalized cycles (large cycles including Hamilton and dominating cycles as special cases) in an arbitrary $k$-connected graph $(k=1,2,...)$ and new lower bounds for the circumference (the length of a longest cycle) are derived, inspiring a number of modifications of famous conjectures of  Bondy (1980) and Jung (2001). All results (both old and new) are shown to be best possible in a sense  based on three types of sharpness  indicating the intervals in which the result is sharp and the intervals in which the result can be further improved. In addition, the presented modifications cannot be derived directly from Bondy's and Jung's conjectures as special cases.\\ 

{\bf MSC-class:} 05C38 (primary), 05C45, 05C40 (secondary)
\end{abstract}

\section{Introduction}
We consider only finite undirected graphs without loops or multiple edges. Notation and terminology not defined here follow that in \cite {3}. Let $G=(V,E)$ be a graph with vertex set $V(G)$ and edge set $E(G).$ For a subset $S$ of $V(G),$ we denote by $G-S$ the maximum subgraph of $G$ with vertex set $V(G)-S.$  For a subgraph $H$ of $G,$ we use $G-H$ to denote $G-V(H).$ 

    Let $\delta$ and $\alpha$  be the  the minimum degree and the independence number of  $G,$ respectively. We define $\sigma_k$ by the  minimum degree sum of any $k$ independent vertices in a graph  if  $\alpha\ge k.$ If $\alpha<k,$ we set $\sigma_k=+\infty.$ In particular, we have $\sigma_1=\delta.$

    By the standard definition, a sequence $v_1 v_2...v_t v_1$ of distinct vertices $v_1,...,v_t$ in $G$ is called a  simple cycle (or just a cycle) of order $t$ (the number of vertices) if $v_i v_{i+1}\in E(G)$ for each $i\in\{1,...,t\},$ where $v_{t+1}=v_1.$ In particular, for $t=1,$ the cycle $v_1$ coincides with the vertex $v_1.$ So, each vertex and edge in a graph can be considered as a cycle of orders 1 and 2, respectively.   

     A cycle in a graph $G$ is called a Hamilton cycle of $G$ if it contains all  the vertices of $G.$  A graph $G$ is called hamiltonian if $G$ contains a Hamilton cycle.

A cycle $Q$ of a graph  $G$ is called dominating cycle  if every edge of $G$ is incident with at least one vertex of $Q.$

In 1980,  Bondy \cite {4} introduced the first type of generalized cycles, including Hamilton and dominating cycles as special cases: for a positive integer $\lambda,$ we call $Q$ a $D_\lambda$-cycle  if $|H|\leq\lambda-1$ for every component $H$ of $G-Q.$ In other words, $Q$ is a $D_\lambda$-cycle of $G$ if and only if every connected subgraph of order $\lambda$ of $G$ has at least one vertex with $Q$ in common. In fact, a $D_\lambda$-cycle dominates all connected subgraphs of order $\lambda.$ According to this definition,  $Q$ is a Hamilton cycle if and only if $Q$ is a $D_1$-cycle; and $Q$ is a dominating cycle if and only if $Q$ is a $D_2$-cycle. 

In this paper, we consider another two types of  generalized cycles  including Hamilton and dominating cycles as special cases. For a positive integer $\lambda$, a cycle $Q$ in $G$ is called a $PD_\lambda$-cycle (PD - Path Dominating) if each path of order at least $\lambda$  has at least one  vertex with $Q$ in common.  Next, we call a cycle $Q$ a $CD_\lambda$-cycle (CD - Cycle Dominating; introduced in  \cite {13}) if each cycle of order at least $\lambda$ has at least one vertex with $Q$ in common. Actually,  a $PD_\lambda$-cycle dominates all paths of order $\lambda$ in $G$; and a $CD_\lambda$-cycle dominates all cycles of order $\lambda$ in $G$. In terms of $PD_\lambda$ and $CD_\lambda$-cycles, $Q$ is a Hamilton cycle if and only if either $Q$ is a $PD_1$-cycle or a $CD_1$-cycle. Further, $Q$ is a dominating cycle if and only if either $Q$ is a $PD_2$-cycle or a $CD_2$-cycle. 

     Throughout the paper, we consider a graph $G$ on $n$ vertices with minimum degree $\delta$ and connectivity $\kappa.$ Further, let $C$ be a longest cycle in $G$ with $c=|C|,$ and let $\overline{p}$ and $\overline{c}$ denote the orders of a longest path and a longest cycle in $G-C,$ respectively. In terms of $\overline{p}$ and $\overline{c},$ $C$ is a Hamilton cycle if and only if either $\overline{p}\leq 0$ or $\overline{c}\leq 0.$ Similarly, $C$ is a dominating cycle if and only if either $\overline{p}\leq 1$ or $\overline{c}\leq 1.$

     In 1980, Bondy \cite {4} conjectured a common generalization of some well-known degree-sum conditions for $PD_\lambda$-cycles (called $(\sigma,\overline{p}$)-version) including Hamilton cycles ($PD_1$-cycles) and dominating cycles ($PD_2$-cycles) as special cases.

\begin{conjecture} (Bondy \cite {4},1980): $(\sigma, \overline{p})$-version. \\
Let $C$ be a longest cycle in a $\lambda$-connected $(1\leq\lambda\leq \delta)$ graph $G$ of order $n.$ If $\sigma_{\lambda+1}\geq n+\lambda(\lambda-1),$ then $\overline{p}\leq \lambda-1.$
\end{conjecture}

    Parts of Conjecture A were proved by Ore \cite {16} $(\lambda=1),$ Bondy \cite {4} $(\lambda=2)$ and Zou \cite {18} $(\lambda=3).$ For the general case, Conjecture A is still open.

     The long cycles analogue (the so called reverse version) of Bondy's conjecture (Conjecture A) can be formulated as follows. 

\begin{conjecture}: (reverse, $\sigma,$ $\overline{p})$-version. \\
Let $C$ be a longest cycle in a $\lambda$-connected $(1\leq\lambda\leq\delta)$ graph $G.$ If $\overline{p}\geq \lambda-1,$ then $c\geq \sigma_\lambda-\lambda(\lambda-2).$
\end{conjecture}

   Parts of Conjecture B were proved by Dirac \cite {6} $(\lambda=1),$ Bondy \cite {2}, Bermond \cite {1}, Linial \cite {11} $(\lambda=2),$ Fraisse, Yung \cite {8} $(\lambda=3)$ and Chiba, Tsugaki, Yamashita \cite {5} $(\lambda=4).$

     The initial motivations of Conjecture A and Conjecture B come from their minimum degree versions - the most popular and much studied versions, which also remain unsolved.  

\begin{conjecture} (Bondy \cite {4},1980): $(\delta, \overline{p})$-version. \\
Let $C$ be a longest cycle in a $\lambda$-connected $(1\leq\lambda\leq \delta)$ graph $G$ of order $n.$ If 
$$
\delta\geq\frac{n+2}{\lambda+1}+\lambda-2,
$$ 
then $\overline{p}\leq \lambda-1.$
\end{conjecture}

Parts of Conjecture C were proved by Dirac \cite {6} $(\lambda=1),$ Nash-Williams \cite{12} $(\lambda=2)$ and  Fan \cite{7} $(\lambda=3).$

\begin{conjecture} (Jung \cite{10}, 2001): (reverse, $\delta$, $\overline{p})$-version. \\
Let $C$ be a longest cycle in a $\lambda$-connected $(1\leq\lambda\leq\delta)$ graph $G.$ If $\overline{p}\geq \lambda-1,$ then $c\geq\lambda(\delta-\lambda+2).$
\end{conjecture}

 Parts of Conjecture D were proved by Dirac \cite{6} $(\lambda=1),$  $(\lambda=2),$ Voss, Zuluaga  \cite{17} $(\lambda=3)$ and Jung \cite {9} $(\lambda=4).$

Note that $CD_\lambda$-cycles are more suitable for research than $PD_\lambda$-cycles since cycles in $G-C$ are more symmetrical than paths in view of the connections between $G-C$ and $CD_\lambda$-cycles. This is the main reason why some minimum degree versions of Conjectures C and D have been solved just for $CD_\lambda$-cycles. 

     According to above arguments, it is natural to consider the exact analogues of Bondy's generalized conjecture (Conjecture A) and its reverse version (Conjecture B) for $CD_\lambda$-cycles which we call $(\sigma,\overline{c})$ and $(reverse, \sigma, \overline{c})$-versions, respectively.

\begin{conjecture}: $(\sigma, \overline{c})$-version. \\
Let $C$ be a longest cycle in a $\lambda$-connected $(1\le\lambda\le\delta)$ graph $G$ of order $n.$ If $\sigma_{\lambda+1}\geq n+\lambda(\lambda-1),$ then $\overline{c}\leq \lambda-1.$
\end{conjecture}

\begin{conjecture}: (reverse, $\sigma$, $\overline{c})$-version. \\
Let $C$ be a longest cycle in a $\lambda$-connected $(1\le\lambda\le\delta)$ graph $G$. If $\overline{c}\geq \lambda-1,$ then $c\geq \sigma_\lambda-\lambda(\lambda-2).$
\end{conjecture}

   In 2009, the author proved  \cite {14} the validity of minimum degree versions of Conjectures E and F.

\begin{theorem} (\cite {14}, 2009): $(\delta, \overline{c})$-version. \\
Let $C$ be a longest cycle in a $\lambda$-connected $(1\le\lambda\le \delta)$ graph $G$ of order $n.$ If 
$$
\delta\geq \frac{n+2}{\lambda+1}+\lambda-2,
$$
then $\overline{c}\leq \lambda-1.$
\end{theorem}

\begin{theorem} (\cite {14}, 2009): (reverse, $\delta$, $\overline{c})$-version. \\
Let $C$ be a longest cycle in a $\lambda$-connected $(1\le\lambda\le\delta)$ graph $G$. If $\overline{c}\geq \lambda-1,$ then $c\geq \lambda(\delta-\lambda+2).$
\end{theorem}

    Actually, in \cite{14} it was proved a significantly stronger result than Theorem A by showing that the conclusion $\overline{c}\le\lambda-1$ in Theorem A can be strengthened to $\overline{c}\le\min\{\lambda-1,\delta-\lambda\},$ called $\overline{c}$-improvement.

\begin{theorem} (\cite {14}, 2009): $(\delta, \overline{c})$-version, $\overline{c}$-improvement. \\
Let $C$ be a longest cycle in a $\lambda$-connected $(1\le\lambda\le\delta)$ graph $G$ of order $n$. If 
$$
\delta\geq \frac{n+2}{\lambda+1}+\lambda-2,
$$ 
then $\overline{c}\leq \min\{ \lambda-1,\delta-\lambda\}.$
\end{theorem}

Recently,  further improvements  of Theorems A and C are presented \cite{15}  inspiring new conjectures in forms of improvements of initial generalized conjecture of Bondy.

\begin{theorem} (\cite{15}, 2022): $(\delta, \overline{c})$-version, $ \kappa$-improvement. \\
Let $C$ be a longest cycle in a graph $G$ of order $n$ and $\lambda$ a positive integer with $1\le\lambda\le\delta.$ If $\kappa\geq \min\{\lambda,\delta-\lambda+1\}$ and 
$$
\delta\geq \frac{n+2}{\lambda+1}+\lambda-2,
$$
then $\overline{c}\leq \lambda-1.$
\end{theorem}

\begin{theorem} (\cite{15}, 2022): $(\delta, \overline{c})$-version, $(\overline{c}, \kappa$)-improvement. \\
Let $C$ be a longest cycle in a graph $G$ of order $n$ and $\lambda$ a positive integer with $1\le\lambda\le\delta.$ If $\kappa\geq \min\{\lambda,\delta-\lambda+1\}$ and 
$$\delta\geq \frac{n+2}{\lambda+1}+\lambda-2,
$$
then $\overline{c}\le\min\{\lambda-1,\delta-\lambda\}.$
\end{theorem}

 Analogously, the condition $\overline{c}\geq \lambda-1$ in Theorem B was weakened  \cite {14} to $\overline{c}\geq \min\{\lambda-1,\delta-\lambda+1\}.$ 

\begin{theorem} (\cite {14}, 2009): (reverse, $\delta$, $\overline{c})$-version, $\overline{c}$-improvement. \\
Let $C$ be a longest cycle in a $\lambda$-connected $(1\le\lambda\le\delta)$ graph $G$. If $\overline{c}\geq\min\{ \lambda-1,\delta-\lambda+1\},$ then $c\geq \lambda(\delta-\lambda+2).$
\end{theorem}

Furthermore,  it was proved \cite{15} that the connectivity condition $\kappa\ge\lambda$ in Theorem B can be weakened to $\kappa\ge\min\{\lambda,\delta-\lambda+2\}.$

\begin{theorem} (\cite{15}, 2022): (reverse, $\delta$, $\overline{c})$-version, $\kappa$-improvement. \\
Let $C$ be a longest cycle in a graph $G$ and $\lambda$ a positive integer with  $1\le\lambda\le\delta.$ If $\kappa\ge\min\{\lambda, \delta-\lambda+2\}$ and $\overline{c}\ge \lambda-1,$ then $c\ge\lambda(\delta-\lambda+2).$  
\end{theorem}

Finally,  it was proved \cite{15} that the connectivity condition $\kappa\ge\lambda$ in Theorem F can be weakened to $\kappa\ge\min\{\lambda,\delta-\lambda+2\}.$

\begin{theorem} (\cite{15}, 2022): (reverse, $\delta$, $\overline{c})$-version, $(\overline{c},\kappa)$-improvement. \\
Let $C$ be a longest cycle in a graph $G$ and $\lambda$ a positive integer with  $1\le\lambda\le\delta.$ If $\kappa\ge\min\{\lambda, \delta-\lambda+2\}$ and $\overline{c}\ge \min\{\lambda-1,\delta-\lambda+1\},$ then $c\ge \lambda(\delta-\lambda+2).$  
\end{theorem}

In this paper, we present a number of new sufficient conditions for large cycles and new lower bounds for the circumference around conjectures of Bondy and Jung which cannot be derived directly from these conjectures as special cases. They can be characterized as  modifications around conjectures of Bondy and Jung.  Each of these modifications is shown to be best possible in a sense. Consider Theorem  B to determine some kinds of sharpness. We say that Theorem B is $c$-sharp if the condition $c\ge \lambda(\delta-\lambda+2)$ cannot be strengthened to $c\ge \lambda(\delta-\lambda+2)+1.$ Further, Theorem B is   $\overline{c}$-sharp if the condition $\overline{c}\ge \lambda-1$ in Theorem B cannot be weakened to $\overline{c}\ge \lambda-2.$  Finally,  Theorem B is $\kappa$-sharp if the  connectivity condition $\kappa\ge \lambda$ cannot be weakened to $\kappa\ge \lambda-1.$ The $\overline{c}$-sharpness,  $\delta$-sharpness and $\kappa$-sharpness  for long cycle versions (reverse versions), say for theorem A, can be defined analogously.

 \section{Modifications}

The first modification can be obtained from Theorem A by replacing the conclusion $\overline{c}\le \lambda-1$ with $\overline{c}\le \delta-\lambda,$ called $\overline{c}$-modification.\\

\noindent {\bf Theorem 1}: $(\delta, \overline{c})$-version, $\overline{c}$-modification. \\
Let $C$ be a longest cycle in a $\lambda$-connected $(1\le\lambda\le\delta)$ graph $G$ of order $n.$ If
 $$
\delta\geq \frac{n+2}{\lambda+1}+\lambda-2,
$$ 
then $\overline{c}\leq \delta-\lambda.$\\

The next  modification can be obtained by replacing the connectivity condition $\kappa\ge \lambda$ in Theorem A with $\kappa\ge \delta-\lambda+1,$ called $ \kappa$-modification.\\

\noindent {\bf Theorem 2}: $(\delta, \overline{c})$-version, $ \kappa$-modification. \\
Let $C$ be a longest cycle in a $(\delta-\lambda+1)$-connected $(1\le\lambda\le\delta)$ graph $G$  of order $n.$ If 
$$
\delta\geq \frac{n+2}{\lambda+1}+\lambda-2,
$$
 then $\overline{c}\leq  \lambda-1.$\\

Another modification can be obtained  by replacing the conclusion $\overline{c}\le \lambda-1$ in Theorem 2 with $\overline{c}\le \delta-\lambda,$  called $ (\overline{c},\kappa$)-modification.\\

\noindent {\bf Theorem 3}: $(\delta, \overline{c})$-version, $ (\overline{c},\kappa$)-modification. \\
Let $C$ be a longest cycle in a $(\delta-\lambda+1)$-connected $(1\le\lambda\le\delta)$ graph $G$ of order $n.$  If  
$$\delta\geq \frac{n+2}{\lambda+1}+\lambda-2,
$$
 then $\overline{c}\leq \delta-\lambda.$\\

Further, Theorem 2 can be improved by strengthening the conclusion $\overline{c}\le \lambda-1$ to $\overline{c}\le \min\{\lambda-1,\delta-\lambda\},$ called $\kappa$-modification with $\overline{c}$-improvement. \\

\noindent{\bf Theorem 4}: $(\delta, \overline{c})$-version,  $\kappa$-modification, $\overline{c}$-improvement. \\
Let $C$ be a longest cycle in a $(\delta-\lambda+1)$-connected $(1\le\lambda\le\delta)$ graph $G$ of order $n.$  If  
$$\delta\geq \frac{n+2}{\lambda+1}+\lambda-2,
$$
 then $\overline{c}\le\min\{\lambda-1,\delta-\lambda\}.$\\

Finally, Theorem 3 can be improved by relaxing the connectivity condition $\kappa\geq \delta-\lambda+1$ to $\kappa\geq \min\{\lambda,\delta-\lambda+1\},$ called $\overline{c}$-modification with $\kappa$-improvement .\\

\noindent{\bf Theorem 5}: $(\delta, \overline{c})$-version, $\overline{c}$-modification, $\kappa$-improvement. \\
Let $C$ be a longest cycle in a graph $G$  of order $n$ and $\lambda$ a positive integer with $1\le\lambda\le\delta.$ If $\kappa\geq \min\{\lambda,\delta-\lambda+1\}$ and 
$$
\delta\geq \frac{n+2}{\lambda+1}+\lambda-2,
$$ 
then $\overline{c}\le\delta-\lambda.$\\

The reverse versions of Theorems 1-5 are generated from Theorem B using different combinations of conditions 
$$
\kappa\ge \lambda,\ \kappa\ge \delta-\lambda+2, \ \kappa\ge\min\{\lambda,\delta-\lambda+2\}
$$
$$
\overline{c}\ge \lambda-1,\ \overline{c}\ge \delta-\lambda+1,\ \overline{c}\ge \min\{\lambda-1,\delta-\lambda+1\}.
$$

The first result can be considered as  $\overline{c}$-modification of Theorem B.\\

\noindent{\bf Theorem 6}: (reverse, $\delta$, $\overline{c})$-version, $\overline{c}$-modification.\\
Let $C$ be a longest cycle in a $\lambda$-connected $(1\le\lambda\le\delta)$ graph G. If  $\overline{c}\ge \delta- \lambda+1,$ then $c\ge\lambda(\delta-\lambda+2).$  \\

The  $\kappa$-modification of Theorem B can be formulated as follows.\\

\noindent{\bf Theorem 7}: (reverse, $\delta$, $\overline{c})$-version, $\kappa$-modification. \\
Let $C$ be a longest cycle in a $(\delta-\lambda+2)$-connected $(1\le\lambda\le\delta)$ graph $G$. If  $\overline{c}\ge \lambda-1,$ then $c\ge\lambda(\delta-\lambda+2).$  \\

Theorem B after $\overline{c}$-modification and $\kappa$-modification, can be formulated as follows.\\

\noindent{\bf Theorem 8}: (reverse, $\delta$, $\overline{c})$-version, $(\overline{c},\kappa)$-modification. \\
Let $C$ be a longest cycle in a $(\delta-\lambda+2)$-connected $(1\le\lambda\le\delta)$ graph $G$. If  $\overline{c}\ge \delta- \lambda+1,$ then $c\ge\lambda(\delta-\lambda+2).$  \\

The next result can be obtained from Theorem 8 by $\overline{c}$-improvement.\\

\noindent{\bf Theorem 9}: (reverse, $\delta$, $\overline{c})$-version, $\kappa$-modification, $\overline{c}$-improvement. \\
Let $C$ be a longest cycle in a $(\delta-\lambda+2)$-connected $(1\le\lambda\le\delta)$ graph $G$. If  $\overline{c}\ge \min\{\lambda-1,\delta-\lambda+1\},$ then $c\ge \lambda(\delta-\lambda+2).$  \\

Finally, the $\kappa$-improvement of Theorem 8 implies the following.\\

\noindent{\bf Theorem 10}: (reverse, $\delta$, $\overline{c})$-version, $\kappa$-improvement, $\overline{c}$-modification. \\
Let $C$ be a longest cycle in a graph $G$ and $\lambda$ a positive integer with  $1\le\lambda\le\delta.$ If $\kappa\ge\min\{\lambda, \delta-\lambda+2\}$ and $\overline{c}\ge \delta-\lambda+1,$ then $c\ge \lambda(\delta-\lambda+2).$\\  

Observe that none of Theorems 1-10 follows from conjectures of Bondy and Jung as a special case.

\section{Generalized modifications}

In this section, motivated by Theorems 1-10 (minimum degree versions), we propose their generalized versions in terms of degree sums ($(\sigma, \overline{c})$-versions)  as generalized modifications around Conjectures of Bondy and Jung.

The first conjecture is a $(\sigma, \overline{c})$-generalization of Theorem 1.\\

\noindent {\bf Conjecture 1}: $(\sigma, \overline{c})$-version, $\overline{c}$-modification. \\
Let $C$ be a longest cycle in a $\lambda$-connected $(1\le\lambda\le\delta)$ graph $G$  of order $n.$ If  $\sigma_{\lambda+1}\geq n+\lambda(\lambda-1),$ then $\overline{c}\leq \delta-\lambda.$\\

The $(\sigma, \overline{c})$-version of Theorem 2 can be formulated as follows.\\

\noindent {\bf Conjecture 2}: $(\sigma, \overline{c})$-version, $\kappa$-modification. \\
Let $C$ be a longest cycle in a $(\delta-\lambda+1)$-connected $(1\le\lambda\le\delta)$ graph $G$ of order $n.$ If  $\sigma_{\lambda+1}\geq n+\lambda(\lambda-1),$ then $\overline{c}\leq \lambda-1.$\\

After $(\sigma, \overline{c})$-generalization, Theorem 3 implies the following.\\

\noindent {\bf Conjecture 3}: $(\sigma, \overline{c})$-version, $\kappa$-modification, $\overline{c}$-modification. \\
Let $C$ be a longest cycle in a $(\delta-\lambda+1)$-connected $(1\le\lambda\le\delta)$ graph $G$ of order $n.$ If  $\sigma_{\lambda+1}\geq n+\lambda(\lambda-1),$ then $\overline{c}\leq \delta-\lambda.$\\

Theorem 4, after $(\sigma, \overline{c})$-generalization, implies the following.\\

\noindent {\bf Conjecture 4}: $(\sigma, \overline{c})$-version, $\kappa$-modification, $\overline{c}$-improvement. \\
Let $C$ be a longest cycle in a $(\delta-\lambda+1)$-connected $(1\le\lambda\le\delta)$ graph $G$ of order $n.$ If  $\sigma_{\lambda+1}\geq n+\lambda(\lambda-1),$ then $\overline{c}\leq \min\{\lambda-1,\delta-\lambda\}.$\\

Finally, Theorem 5 implies the next conjecture after $(\sigma, \overline{c})$-generalization.\\

\noindent {\bf Conjecture 5}: $(\sigma, \overline{c})$-version, $\overline{c}$-modification, $\kappa$-improvement. \\
Let $C$ be a longest cycle in a $\min\{\lambda,\delta-\lambda+1\}$-connected $(1\le\lambda\le\delta)$ graph $G$  of order $n.$ If  $\sigma_{\lambda+1}\geq n+\lambda(\lambda-1),$ then $\overline{c}\leq \delta-\lambda.$\\

The next five conjectures can be proposed by replacing the conclusion $c\ge \lambda(\delta-\lambda+2)$ in Theorems 6-10 with $c\ge\sigma_\lambda-\lambda(\lambda-2).$\\

\noindent{\bf Conjecture 6}: (reverse, $\sigma$, $\overline{c})$-version, $\overline{c}$-modification.\\
Let $C$ be a longest cycle in a $\lambda$-connected $(1\le\lambda\le\delta)$ graph $G$. If  $\overline{c}\ge \delta- \lambda+1,$ then $c\ge\sigma_\lambda-\lambda(\lambda-2).$  \\

\noindent{\bf Conjecture 7}: (reverse, $\sigma$, $\overline{c})$-version, $\kappa$-modification.\\
Let $C$ be a longest cycle in a $(\delta-\lambda+2)$-connected $(1\le\lambda\le\delta)$ graph $G$. If  $\overline{c}\ge  \lambda-1,$ then $c\ge\sigma_\lambda-\lambda(\lambda-2).$  \\

\noindent{\bf Conjecture 8}: (reverse, $\sigma$, $\overline{c})$-version, $\kappa$-modification, $\overline{c}$-modification.\\
Let $C$ be a longest cycle in a $(\delta-\lambda+2)$-connected $(1\le\lambda\le\delta)$ graph $G$. If  $\overline{c}\ge \delta- \lambda+1,$ then $c\ge\sigma_\lambda-\lambda(\lambda-2).$  \\

\noindent{\bf Conjecture 9}: (reverse, $\sigma$, $\overline{c})$-version, $\kappa$-modification, $\overline{c}$-improvement.\\
Let $C$ be a longest cycle in a $(\delta-\lambda+2)$-connected $(1\le\lambda\le\delta)$ graph $G$. If  $\overline{c}\ge \min\{\lambda-1,\delta- \lambda+1\},$ then $c\ge\sigma_\lambda-\lambda(\lambda-2).$  \\

\noindent{\bf Conjecture 10}: (reverse, $\sigma$, $\overline{c})$-version, $\overline{c}$-modification, $\kappa$-improvement.\\
Let $C$ be a longest cycle in a $\min\{\lambda,\delta-\lambda+2\}$-connected $(1\le\lambda\le\delta)$ graph $G$. If  $\overline{c}\ge \delta- \lambda+1,$ then $c\ge\sigma_\lambda-\lambda(\lambda-2).$  \\

The $(\sigma, \overline{p})$-versions of Conjectures 1-10 can be proposed by a similar way.

 \section{Proofs}

The proofs are based on several implications with the following two schemes:․\\

$$
\text{Theorem E \ \  $\Rightarrow$} \ \  \left\{
\begin{array}{lll}
\text{Theorem C \ \  $\Rightarrow$ \ \ Theorem A,} \\
\text{Theorem D \ \ $\Rightarrow$ \ \ Theorem A,} \\
\text{Theorem 4 \ \ \  $\Rightarrow$ \ \  Theorem 2 and Theorem 3}. \\
\text{Theorem 5 \ \ \  $\Rightarrow$ \ \  Theorem 1 and Theorem 3,} \\
\end{array}
\right.
$$

$$
\text{Theorem H \ \  $\Rightarrow$} \ \  \left\{
\begin{array}{lll}
\text{Theorem F \ \ $\Rightarrow$ \ \ Theorem B,} \\
\text{Theorem G \ \  $\Rightarrow$ \ \ Theorem B,} \\
\text{Theorem 9 \ \ \  $\Rightarrow$ \ \  Theorem 7 and Theorem 8,} \\
\text{Theorem 10 \ \  $\Rightarrow$ \ \  Theorem 6 and Theorem 8}. \\
\end{array}
\right.
$$

 \section{On sharpness}

Let $m,t$ be positive integers. Take $t+1$ disjoint copies of the complete graph $K_m$ and join each vertex in their union to every vertex of a disjoint complete graph $K_t,$ denoted by $G=(t+1)K_m+K_t.$

To show the $c,\overline{c},\kappa,\delta$-sharpness of all modifications (Theorems 1-10), we need analogous observations around  original Theorems A, B and for all improvements: Theorems C, D, E, F, G and H. We will use different graph examples of standard form $G=(t+1)K_m+K_t$ for appropriate $m$ and $t$ guaranteeing the minimum degree $\delta(G)=m+t-1.$ 

In particular, for the $\overline{c}$-sharpness of Theorem A, take   $G=(\delta-\lambda+3)K_{\lambda-1}+K_{\delta-\lambda+2}$ with $n(G)=\lambda(\delta-\lambda+3)-1.$  The connectivity condition $\kappa(G)=\delta-\lambda+2\ge\lambda$ and the minimum degree condition
$$
\delta\ge \frac{n+2}{\lambda+1}+\lambda-2    \eqno{(1)}
$$
follow immediately whenever $\lambda\le\frac{\delta+2}{2}$. Observing also that $\overline{c}(G)=\lambda-1,$ we conclude that Theorem A is $\overline{c}$-sharp for each  $\lambda\le\frac{\delta+1}{2}$.

For the $\kappa$-sharpness of Theorem A, take $G=\lambda K_{\delta-\lambda+2}+K_{\lambda-1}$ with $n(G)=\lambda(\delta-\lambda+3)-1.$  If $\lambda\le\frac{\delta+1}{2}$, then the condition (1) and the inequality  $\overline{c}(G)=\delta-\lambda+2\ge \lambda$ follow immediately. Since $\kappa(G)=\lambda-1,$ we conclude that the connectivity condition $\kappa\ge \lambda$ in Theorem A cannot be weakened to $\kappa\ge \lambda-1$. In other words, Theorem A is $\kappa$-sharp for each $\lambda\le\frac{\delta+1}{2}.$  

For the  $\delta$-sharpness of Theorem A, take $G=(\delta-\lambda+2)K_{\lambda}+K_{\delta-\lambda+1}$  on $n(G)=(\lambda+1)(\delta-\lambda+2)-1$ vertices. Then
$$
\delta=\frac{(\lambda+1)(\delta-\lambda+2)}{\lambda+1}+\lambda-2=\frac{n+1}{\lambda+1}+\lambda-2.     \eqno{(2)}
$$ 

Since the connectivity condition $\kappa(G)=\delta-\lambda+1\ge\lambda$ holds for each $\lambda\le\frac{\delta+1}{2}$ and  $\overline{c}(G)=\lambda>\lambda-1,$ we conclude that Theorem A is $\delta$-sharp for each $\lambda\le\frac{\delta+1}{2}.$

Combining above observations, we obtain the following.\\

 \noindent {\bf Proposition 1}. 
$$
\text{Theorem A is} \ \  \left\{
\begin{array}{lll}
\text{$\overline{c}$-sharp for each $\lambda\leq \frac{\delta +1}{2},$} \\
\text{$\kappa$-sharp for each $\lambda\leq \frac{\delta +1}{2},$} \\
\text{$\delta$-sharp for each $\lambda\leq \frac{\delta +1}{2}.$} \\
\end{array}
\right.
$$\\

Now consider the sharpness of Theorem 3. For the $\overline{c}$-sharpness, take $G=(\lambda+2)K_{\delta-\lambda}+K_{\lambda+1}.$ Then   $n(G)=(\lambda+2)(\delta-\lambda+1)-1.$ It is easy to see that the connectivity condition $\kappa(G)=\lambda+1\ge\delta-\lambda+1$ and the minimum degree condition (1) are guaranteed whenever $\lambda\ge\frac{\delta+1}{2}.$ Since $\overline{c}(G)=\delta-\lambda,$ Theorem 3 is $\overline{c}$-sharp for each $\lambda\ge\frac{\delta+1}{2}.$

For the $\kappa$-sharpness, take $G=(\delta-\lambda+1)K_{\lambda+1}+K_{\delta-\lambda}$ on  $n(G)=(\lambda+2)(\delta-\lambda+1)-1$ vertices.  The inequality $\overline{c}(G)=\lambda+1\ge \delta-\lambda+1$ and the minimum degree condition (1) are guaranteed  whenever $\lambda\ge\frac{\delta+1}{2}.$ Observing that $\kappa(G)=\delta-\lambda,$ we conclude that Theorem 3 is $\kappa$-sharp for each  $\lambda\ge\frac{\delta+1}{2}.$

For the $\delta$-sharpness, take $G=(\delta-\lambda+2)K_{\lambda}+K_{\delta-\lambda+1}.$ Then  
$$
\kappa(G)=\delta-\lambda+1, \ n(G)=(\lambda+1)(\delta-\lambda+2)-1.
$$

Therefore, the equality (2) follows immediately. Observing also that the inequality $\overline{c}(G)=\lambda\ge \delta-\lambda+1$ holds whenever $\lambda\ge\frac{\delta+1}{2},$ we conclude that Theorem 3 is $\delta$-sharp for each $\lambda\ge\frac{\delta+1}{2}.$ 

Combining above observations, we can claim  the following.\\

\noindent {\bf Proposition 2}. 
$$
\text{Theorem 3 is} \ \  \left\{
\begin{array}{lll}
\text{$\overline{c}$-sharp for each $\lambda\geq \frac{\delta +1}{2},$} \\
\text{$\kappa$-sharp for each $\lambda\geq \frac{\delta +1}{2},$} \\
\text{$\delta$-sharp for each $\lambda\geq \frac{\delta +1}{2}.$} \\
\end{array}
\right.
$$\\

Consider the sharpness of Theorem 1. For the $\overline{c}$-sharpness, take $G=(\lambda+2)K_{\delta-\lambda}+K_{\lambda+1}.$ Then 
$$
 \kappa(G)=\lambda+1>\lambda, \ n(G)=(\lambda+2)(\delta-\lambda+1)-1.
$$

It is easy to see that the  minimum degree condition (1) is guaranteed whenever $\lambda\ge\frac{\delta+1}{2}.$ Since $\overline{c}(G)=\delta-\lambda,$ Theorem 1 is $\overline{c}$-sharp for each $\lambda\ge\frac{\delta+1}{2}.$

For the $\kappa$-sharpness, take $G=\lambda K_{\delta-\lambda+2}+K_{\lambda-1}$ on $n(G)=\lambda(\delta-\lambda+3)-1$ vertices. If $\lambda\le\frac{\delta+1}{2},$ then the minimum degree condition (1) holds immediately. Furthermore, since $\kappa(G)=\lambda-1$ and 
$$
\overline{c}(G)=\delta-\lambda+2>\delta-\lambda+1,
$$
Theorem 1  is $\kappa$-sharp for each $\lambda\le\frac{\delta+1}{2}.$  

For the $\delta$-sharpness, take  $G=(\lambda+1)K_{\delta-\lambda+1}+K_{\lambda}.$ Then 
$$
\kappa(G)=\lambda, \ n(G)=(\lambda+1)(\delta-\lambda+2)-1.
$$

The equality (2) follows immediately. Observing also that $\overline{c}(G)=\delta-\lambda+1,$ we conclude that Theorem 1 is $\delta$-sharp for each $1\le\lambda\le\delta.$ 

Thus, we have the following. \\

\noindent {\bf Proposition 3}.

$$
\text{Theorem 1 is} \ \  \left\{
\begin{array}{lll}
\text{$\overline{c}$-sharp for each $\lambda\geq \frac{\delta +1}{2},$} \\
\text{$\kappa$-sharp for each $\lambda\leq \frac{\delta +1}{2},$} \\
\text{$\delta$-sharp for each $1\le\lambda\le\delta.$} \\
\end{array}
\right.
$$\\

For Theorem 2,  take $G=(\delta-\lambda+3)K_{\lambda-1}+K_{\delta-\lambda+2}.$ Then 
$$
n(G)=\lambda(\delta-\lambda+3)-1, \ \kappa(G)=\delta-\lambda+2>\delta-\lambda+1.
$$

The  minimum degree condition (1)  is guaranteed whenever $\lambda\le\frac{\delta+1}{2}.$ Recalling also that $\overline{c}(G)=\lambda-1,$ we conclude that  Theorem 2 is $\overline{c}$-sharp for each  $\lambda\le\frac{\delta+1}{2}.$

For the $\kappa$-sharpness, let $G=(\delta-\lambda+1)K_{\lambda+1}+K_{\delta-\lambda}.$ Then $n(G)=(\lambda+2)(\delta-\lambda+1)-1.$   The  minimum degree condition (1) holds for each  $\lambda\ge\frac{\delta+1}{2}.$ Observing also that $\kappa(G)=\delta-\lambda$ and $\overline{c}(G)=\lambda+1>\lambda,$ we conclude that Theorem 2 is $\kappa$-sharp for each  $\lambda\ge\frac{\delta+1}{2}.$

For the $\delta$-sharpness, let $G=(\delta-\lambda+2)K_{\lambda}+K_{\delta-\lambda+1}.$ Then 
$$
\kappa(G)=\delta-\lambda+1, \   n(G)=(\lambda+1)(\delta-\lambda+2)-1.
$$

The equality (2) follows immediately. Since $\overline{c}(G)=\lambda,$   Theorem 2 is $\delta$-sharp for each $1\le\lambda\le\delta.$ 

Hence, we get the following.\\

\noindent {\bf Proposition 4}. 
$$
\text{Theorem 2 is} \ \  \left\{
\begin{array}{lll}
\text{$\overline{c}$-sharp for each $\lambda\leq \frac{\delta +1}{2},$} \\
\text{$\kappa$-sharp for each $\lambda\geq \frac{\delta +1}{2},$} \\
\text{$\delta$-sharp for each $1\le\lambda\le\delta.$} \\
\end{array}
\right.
$$\\

Consider the sharpness of Theorem C. If $\lambda\le \frac{\delta+1}{2},$ then we have Theorem A. By Proposition 1, Theorem C is $(\overline{c},\kappa,\delta)$-sharp. If $\lambda\ge \frac{\delta+1}{2},$ then we have Theorem 1. By Proposition 3, Theorem C is $(\overline{c},\delta)$-sharp. So, we have the following.\\

\noindent {\bf Proposition 5}. 
$$
\text{Theorem C is} \ \  \left\{
\begin{array}{lll}
\text{$\overline{c}$-sharp for each $1\le\lambda\le\delta,$} \\
\text{$\kappa$-sharp for each $\lambda\leq \frac{\delta +1}{2},$} \\
\text{$\delta$-sharp for each $1\le\lambda\le\delta.$} \\
\end{array}
\right.
$$\\

Consider the sharpness of Theorem D. If $\lambda\le \frac{\delta+1}{2},$ then we have Theorem A. By Proposition 1, Theorem D is $(\overline{c},\kappa,\delta)$-sharp. If $\lambda\ge \frac{\delta+1}{2},$ then we have Theorem 2. By Proposition 4, Theorem D is $(\kappa,\delta)$-sharp. So, we have the following.\\

\noindent {\bf Proposition 6}. 
$$
\text{Theorem D is} \ \  \left\{
\begin{array}{lll}
\text{$\overline{c}$-sharp for each $\lambda\leq \frac{\delta +1}{2},$} \\
\text{$\kappa$-sharp for each $1\le\lambda\le\delta,$} \\
\text{$\delta$-sharp for each $1\le\lambda\le\delta.$} \\
\end{array}
\right.
$$\\

Consider the sharpness of Theorem 5. If $\lambda\le \frac{\delta+1}{2},$ then we have Theorem 1. By Proposition 3, Theorem 5 is $(\kappa,\delta)$-sharp. If $\lambda\ge \frac{\delta+1}{2},$ then we have Theorem 3. By Proposition 2, Theorem 5 is $(\overline{c},\kappa,\delta)$-sharp. So, we have the following.\\

\noindent {\bf Proposition 7}. 

$$
\text{Theorem 5 is} \ \  \left\{
\begin{array}{lll}
\text{$\overline{c}$-sharp for each $\lambda\geq \frac{\delta +1}{2},$} \\
\text{$\kappa$-sharp for each $1\le\lambda\le\delta,$} \\
\text{$\delta$-sharp for each $1\le\lambda\le\delta.$} \\
\end{array}
\right.
$$\\

Consider the sharpness of Theorem 4. If $\lambda\le \frac{\delta+1}{2},$ then we have Theorem 2. By Proposition 4, Theorem 4 is $(\overline{c},\delta)$-sharp. If $\lambda\ge \frac{\delta+1}{2},$ then we have Theorem 3. By Proposition 2, Theorem 4 is $(\overline{c},\kappa,\delta)$-sharp. So, we have the following.\\

\noindent {\bf Proposition 8}. 

$$
\text{Theorem 4 is} \ \  \left\{
\begin{array}{lll}
\text{$\overline{c}$-sharp for each $1\le\lambda\le\delta,$} \\
\text{$\kappa$-sharp for each $\lambda\geq \frac{\delta +1}{2},$} \\
\text{$\delta$-sharp for each $1\le\lambda\le\delta.$} \\
\end{array}
\right.
$$\\

Consider the sharpness of Theorem E. If $\lambda\le \frac{\delta+1}{2},$ then we have Theorem A. By Proposition 1, Theorem E is $(\overline{c},\kappa,\delta)$-sharp. If $\lambda\ge \frac{\delta+1}{2},$ then we have Theorem 3. By Proposition 2, Theorem E is $(\overline{c},\kappa,\delta)$-sharp as well. So, we have the following.\\

\noindent {\bf Proposition 9}. 
$$
\text{Theorem E is} \ \  \left\{
\begin{array}{lll}
\text{$\overline{c}$-sharp for each $1\le\lambda\le\delta,$} \\
\text{$\kappa$-sharp for each $1\le\lambda\le\delta$,} \\
\text{$\delta$-sharp for each $1\le\lambda\le\delta$.} \\
\end{array}
\right.
$$\\

Now we turn to the sharpness observations around Jung's conjecture. 

Consider the sharpness of Theorem B. For the $c$-sharpness, let $G=(\delta-\lambda+3)K_{\lambda-1}+K_{\delta-\lambda+2}.$ Then $\overline{c}(G)=\lambda-1.$   The connectivity condition $\kappa(G)=\delta-\lambda+2\ge\lambda$ holds whenever $\lambda\le \frac{\delta+2}{2}.$ Observing also that $c(G)=\lambda(\delta-\lambda+2),$ we conclude that Theorem B is $c$-sharp for each $\lambda\le \frac{\delta+2}{2}.$

For the $\kappa$-sharpness, take $G=\lambda K_{\delta-\lambda+2}+K_{\lambda-1}$ with $\kappa(G)=\lambda-1$ and  $c(G)=(\lambda-1)(\delta-\lambda+3).$ Since the inequalities 
$$
c(G)<\lambda(\delta-\lambda+2), \ \overline{c}(G)=\delta-\lambda+2\ge\lambda-1 
$$
are satisfied for each $\lambda\le \frac{\delta+2}{2},$  we conclude that Theorem B is $\kappa$-sharp for each  $\lambda\le \frac{\delta+2}{2}.$

For the $\overline{c}$-sharpness, take $G=(\delta-\lambda+4)K_{\lambda-2}+K_{\delta-\lambda+3}$ with  $\overline{c}(G)=\lambda-2$ and $c(G)=(\lambda-1)(\delta-\lambda+3).$ Observing that the  inequalities 
$$
c(G)<\lambda(\delta-\lambda+2), \ \kappa(G)=\delta-\lambda+3\ge\lambda
$$
are satisfied under the condition $\lambda\le \frac{\delta+2}{2},$   we conclude that Theorem B is $\overline{c}$-sharp for each $\lambda\le \frac{\delta+2}{2}.$

Combining these observations, we have the following.\\

\noindent {\bf Proposition 10}. 

$$
\text{Theorem B is} \ \  \left\{
\begin{array}{lll}
\text{$c$-sharp for each $\lambda\le\frac{\delta+2}{2},$} \\
\text{$\kappa$-sharp for each $\lambda\le\frac{\delta+2}{2},$} \\
\text{$\overline{c}$-sharp for each $\lambda\le\frac{\delta+2}{2}.$} \\
\end{array}
\right.
$$\\

For the $c$-sharpness of Theorem 8, take $G=(\delta-\lambda+3)K_{\lambda-1}+K_{\delta-\lambda+2}$ with 
$$
\kappa(G)=\delta-\lambda+2, \ \overline{c}(G)=\lambda-1.
$$
The inequality $\overline{c}(G)=\lambda-1\ge\delta-\lambda+1$ holds whenever $\lambda\ge \frac{\delta+2}{2}.$ Observing also that $c(G)=\lambda(\delta-\lambda+2),$ we conclude that  Theorem 8 is $c$-sharp for each $\lambda\ge \frac{\delta+2}{2}.$

For the $\kappa$-sharpness, let $G=(\delta-\lambda+2)K_{\lambda}+K_{\delta-\lambda+1}$ with 
$$
\kappa(G)=\delta-\lambda+1, \ c(G)=(\lambda+1)(\delta-\lambda+1).
$$
Observing also that the inequalities
$$
c(G)<\lambda(\delta-\lambda+2), \ \overline{c}(G)=\lambda\ge\delta-\lambda+1
$$  
are satisfied whenever $\lambda\ge \frac{\delta+2}{2},$ we conclude that Theorem 8  is $\kappa$-sharp for each $\lambda\ge \frac{\delta+2}{2}.$

For the $\overline{c}$-sharpness, take $G=(\lambda+2)K_{\delta-\lambda}+K_{\lambda+1}$ with 
$$
\overline{c}(G)=\delta-\lambda, \ c(G)=(\lambda+1)(\delta-\lambda+1).
$$
Since the inequalities 
$$
\kappa(G)=\lambda+1\ge\delta-\lambda+2, \ c(G)<\lambda(\delta-\lambda+2)
$$
are satisfied whenever $\lambda\ge\frac{\delta+2}{2},$ we conclude that Theorem 8 is $\overline{c}$-sharp for each  $\lambda\ge\frac{\delta+2}{2}.$

Combining these observations, we have the following.\\

\noindent {\bf Proposition 11}. 
$$
\text{Theorem 8 is} \ \  \left\{
\begin{array}{lll}
\text{$c$-sharp for each $\lambda\ge\frac{\delta+2}{2},$} \\
\text{$\kappa$-sharp for each $\lambda\ge\frac{\delta+2}{2},$} \\
\text{$\overline{c}$-sharp for each $\lambda\ge\frac{\delta+2}{2}.$} \\
\end{array}
\right.
$$\\

For the $c$-sharpness of Theorem 6, take $G=(\lambda+1)K_{\delta-\lambda+1}+K_{\lambda}.$ Then  
$$
\overline{c}(G)=\delta-\lambda+1, \ \kappa(G)=\lambda.
$$

Since $c(G)=\lambda(\delta-\lambda+2),$ Theorem 6 is $c$-sharp for each $1\le\lambda\le\delta.$

For the $\kappa$-sharpness, take $G=\lambda K_{\delta-\lambda+2}+K_{\lambda-1}$ with $\kappa(G)=\lambda-1.$  Then  
$$
\overline{c}(G)=\delta-\lambda+2>\delta-\lambda+1, \ c(G)=(\lambda-1)(\delta-\lambda+3).
$$

Since the inequality $c(G)<\lambda(\delta-\lambda+2)$ is satisfied whenever $\lambda\le\frac{\delta+2}{2},$    Theorem 6 is $\kappa$-sharp for each  $\lambda\le \frac{\delta+2}{2}.$

 For the $\overline{c}$-sharpness, let $G=(\lambda+2)K_{\delta-\lambda}+K_{\lambda+1}$ with $\overline{c}(G)=\delta-\lambda.$  Then
$$
\kappa(G)=\lambda+1>\lambda, \ c(G)=(\lambda+1)(\delta-\lambda+1).
$$
Since the inequality  $c(G)<\lambda(\delta-\lambda+2)$ is satisfied whenever $\lambda\ge\frac{\delta+2}{2},$   Theorem 6 is $\overline{c}$-sharp for each  $\lambda\ge\frac{\delta+2}{2}.$

Combining these observations, we have the following.\\

\noindent {\bf Proposition 12}. 
$$
\text{Theorem 6 is} \ \  \left\{
\begin{array}{lll}
\text{$c$-sharp for each $1\le\lambda\le\delta,$} \\
\text{$\kappa$-sharp for each $\lambda\le\frac{\delta+2}{2},$} \\
\text{$\overline{c}$-sharp for each $\lambda\ge\frac{\delta+2}{2},$} \\
\end{array}
\right.
$$\\

For the $c$-sharpness of Theorem 7, take $G=(\delta-\lambda+3)K_{\lambda-1}+K_{\delta-\lambda+2}$ with  
$$
\overline{c}(G)=\lambda-1, \  \kappa(G)=\delta-\lambda+2.
$$

Since   $c(G)=\lambda(\delta-\lambda+2)$,  Theorem 7 is $c$-sharp for each $1\le\lambda\le \delta.$

For the $\kappa$-sharpness, take  $G=(\delta-\lambda+2)K_{\lambda}+K_{\delta-\lambda+1}$  with 
$$
\kappa(G)=\delta-\lambda+1, \ \overline{c}(G)=\lambda>\lambda-1, \ c(G)=(\lambda+1)(\delta-\lambda+1).
$$
Since the inequality $c(G)<\lambda(\delta-\lambda+2)$ holds whenever  $\lambda\ge\frac{\delta+2}{2},$  Theorem 7 is $\kappa$-sharp for each $\lambda\ge\frac{\delta+2}{2}.$

For the  $\overline{c}$-sharpness, let $G=(\delta-\lambda+4)K_{\lambda-2}+K_{\delta-\lambda+3}.$ Then  
$$
\overline{c}(G)=\lambda-2, \ \kappa(G)=\delta-\lambda+3>\delta-\lambda+2, \ c(G)=(\lambda-1)(\delta-\lambda+3).
$$

Sice the  inequality $c(G)<\lambda(\delta-\lambda+2)$ is satisfied whenever $\lambda\le \frac{\delta+2}{2},$  Theorem 7 is $\overline{c}$-sharp for each $\lambda\le \frac{\delta+2}{2}.$ 

Combining these observations, we get the following.\\

\noindent {\bf Proposition 13}. 
$$
\text{Theorem 7 is} \ \  \left\{
\begin{array}{lll}
\text{$c$-sharp for each $1\le\lambda\le\delta,$} \\
\text{$\kappa$-sharp for each $\lambda\ge\frac{\delta+2}{2},$} \\
\text{$\overline{c}$-sharp for each $\lambda\le\frac{\delta+2}{2}.$} \\
\end{array}
\right.
$$\\

Consider the sharpness of Theorem F. If $\lambda\le \frac{\delta+2}{2},$ then we have Theorem B. By Proposition 10, Theorem F is $(c,\kappa,\overline{c})$-sharp for each $\lambda\le \frac{\delta+2}{2}$. If $\lambda\ge \frac{\delta+2}{2}$, then we have Theorem 6. By Proposition 12, Theorem F is $(c,\overline{c})$-sharp for each $\lambda\ge \frac{\delta+2}{2}$. Combining these observations, we obtain the following.\\ 

\noindent {\bf Proposition 14}. 
$$
\text{Theorem F is} \ \  \left\{
\begin{array}{lll}
\text{$c$-sharp for each $1\le\lambda\le\delta,$} \\
\text{$\kappa$-sharp for each $\lambda\le\frac{\delta+2}{2},$} \\
\text{$\overline{c}$-sharp for each $1\le\lambda\le\delta.$} \\
\end{array}
\right.
$$\\

Consider the sharpness of Theorem G. If $\lambda\le \frac{\delta+2}{2},$ then we have Theorem B. By Proposition 10, Theorem G is $(c,\kappa,\overline{c})$-sharp for each $\lambda\le \frac{\delta+2}{2}.$ If $\lambda\ge \frac{\delta+2}{2},$ then we have Theorem 7. By Proposition 13, Theorem G is $(c,\kappa)$-sharp for each $\lambda\ge \frac{\delta+2}{2}.$ Combining these observations, we obtain the following.\\ 

\noindent {\bf Proposition 15}. 
$$
\text{Theorem G is} \ \  \left\{
\begin{array}{lll}
\text{$c$-sharp for each $1\le\lambda\le\delta,$} \\
\text{$\kappa$-sharp for each $1\le\lambda\le\delta,$} \\
\text{$\overline{c}$-sharp for each $\lambda\le\frac{\delta+2}{2}.$} \\
\end{array}
\right.
$$\\

Consider the sharpness of Theorem 10. If $\lambda\le \frac{\delta+2}{2},$ then we have Theorem 6. By Proposition 12, Theorem 10 is $(c,\kappa)$-sharp for each $\lambda\le \frac{\delta+2}{2}.$ If $\lambda\ge \frac{\delta+2}{2},$ then we have Theorem 8. By Proposition 11, Theorem 10 is $(c,\kappa,\overline{c})$-sharp for each $\lambda\ge \frac{\delta+2}{2}.$ Combining these observations, we obtain the following.\\ 

\noindent {\bf Proposition 16}. 
$$
\text{Theorem 10 is} \ \  \left\{
\begin{array}{lll}
\text{$c$-sharp for each $1\le\lambda\le\delta,$} \\
\text{$\kappa$-sharp for each $1\le\lambda\le\delta,$} \\
\text{$\overline{c}$-sharp for each $\lambda\ge\frac{\delta+2}{2}.$} \\
\end{array}
\right.
$$\\

Consider the sharpness of Theorem 9. If $\lambda\le \frac{\delta+2}{2},$ then we have Theorem 7. By Proposition 13, Theorem 9 is $(c,\overline{c})$-sharp for each $\lambda\le \frac{\delta+2}{2}.$ If $\lambda\ge \frac{\delta+2}{2},$ then we have Theorem 8. By Proposition 11, Theorem 9 is $(c,\kappa,\overline{c})$-sharp for each $\lambda\ge \frac{\delta+2}{2}.$ Combining these observations, we obtain the following.\\

\noindent {\bf Proposition 17}. 
$$
\text{Theorem 9 is} \ \  \left\{
\begin{array}{lll}
\text{$c$-sharp for each $1\le\lambda\le\delta,$} \\
\text{$\kappa$-sharp for each $\lambda\ge\frac{\delta+2}{2},$} \\
\text{$\overline{c}$-sharp for each $1\le\lambda\le\delta.$} \\
\end{array}
\right.
$$\\

Consider the sharpness of Theorem H. If $\lambda\le \frac{\delta+2}{2},$ then we have Theorem B. By Proposition 10, Theorem H is $(c,\kappa,\overline{c})$-sharp for each $\lambda\le \frac{\delta+2}{2}.$ If $\lambda\ge \frac{\delta+2}{2},$ then we have Theorem 8. By Proposition 11, Theorem H is $(c,\kappa,\overline{c})$-sharp for each $\lambda\ge \frac{\delta+2}{2}.$ Combining these observations, we obtain the following.\\ 

\noindent {\bf Proposition 18}. 
$$
\text{Theorem H is} \ \  \left\{
\begin{array}{lll}
\text{$c$-sharp for each $1\le\lambda\le\delta,$} \\
\text{$\kappa$-sharp for each $1\le\lambda\le\delta,$} \\
\text{$\overline{c}$-sharp for each $1\le\lambda\le\delta.$} \\
\end{array}
\right.
$$\\

\noindent  Institute for Informatics and Automation Problems of NAS RA\\
P. Sevak 1, Yerevan 0014, Armenia\\
e-mail: zhora@iiap.sci.am

\end{document}